\documentclass{amsart}
\usepackage[utf8]{inputenc}

\usepackage{tikz}
\usepackage{xcolor}

\usepackage{amsmath}
\usepackage{amsthm}
\usepackage{amssymb}

\usepackage{booktabs}

\usepackage[bookmarks=false]{hyperref}

\newtheorem{thm}{Theorem}[section]

\newtheorem*{ack}{Acknowledgements}

\theoremstyle{definition}

\newtheorem{rem}[thm]{Remark}

\newcommand{\bbC}{\mathbb{C}}
\newcommand{\bbR}{\mathbb{R}}
\newcommand{\bbQ}{\mathbb{Q}}
\newcommand{\bbZ}{\mathbb{Z}}
\newcommand{\bbS}{\mathbb{S}}

\newcommand{\bbT}{\mathbb{T}}
\newcommand{\bbP}{\mathbb{P}}
\newcommand{\bbG}{\mathbb{G}}

\DeclareMathOperator{\Hom}{Hom}
\DeclareMathOperator{\Aut}{Aut}

\DeclareMathOperator{\SO}{SO}
\DeclareMathOperator{\SL}{SL}
\DeclareMathOperator{\Ric}{Ric}

\def\indicator#1{{\mathbf{1}_{#1}}}

\title{Weight sensitivity in K-stability of Fano varieties
}
\author{Thibaut Delcroix}
\address{Thibaut Delcroix, Univ Montpellier, CNRS, Montpellier, France}
\email{thibaut.delcroix@umontpellier.fr}
\date{\today}

\begin{document}

\begin{abstract}
We prove that, for a spherical Fano threefold not in the Mori-Mukai family 2-29, and a weight function associated with the action of the connected center of a Levi subgroup of its automorphism group, weighted K-polystability is equivalent to vanishing of the weighted Futaki invariant.  
This is surprising since unlike the case of toric Fano manifold, there exist non-product, special, equivariant test configurations. 
For the Kähler-Einstein Fano threefold 2-29, and for well-chosen torus action on the three dimensional quadric, we show that this property is false and exhibit explicit examples of weighted optimal degenerations. 
We then generalize this to higher-dimensional quadrics and blowups of quadrics along a codimension 2 subquadric. 
\end{abstract}

\maketitle

\section{Introduction}

K-stability is the algebro-geometric notion corresponding, via the Yau-Tian-Donaldson conjecture, to existence of canonical Kähler metrics on complex varieties. 
In the case of Fano manifolds, the notion of K-stability corresponding to Kähler-Einstein metrics has been tremendously developed following the resolution of the Yau-Tian-Donaldson conjecture in this setting \cite{CDSI,CDSII,CDSIII,Tian_2015}, yielding numerous new examples of manifolds with or without Kähler-Einstein metrics (see e.g. \cite{Ilten-Suss_2017,Delcroix_2020,the_book}), and applications to the study of moduli of Fano varieties (see e.g. Xu's survey \cite{Xu_2021}). 
The references given in the previous sentence are only intended as illustrating landmarks, and do not reflect many of the recent achievements, the subject being extremely active. 
Since many Fano manifolds do not admit Kähler-Einstein metrics, it is desirable to have an alternative notion of canonical Kähler metric for those. 
Kähler-Ricci solitons, obtained as limits of the Ricci flow on Kähler manifolds, provide natural candidates and the study of moduli of Fano manifolds with Kähler-Ricci solitons, or with the corresponding notion of weighted K-polystability has been highlighted by Chenyang Xu as a possible future direction for the study of K-moduli. 

More generally, one can consider, on a Fano manifold \(X\) equipped with a holomorphic action of a compact torus \(T_c\simeq (\bbS^1)^k\), the notion of weighted soliton. 
To define this, consider the canonically normalized moment image \(\Delta\) of the Kähler class \(c_1(X)\) with respect to the \(T_c\)-action. 
Given a smooth, positive function \(g:\Delta \to \bbR_{>0}\), a smooth Kähler form \(\omega\in c_1(X)\) is called a \emph{weighted soliton} with respect to the weight \(g\) (or shortly, a \(g\)-soliton) if 
\[ \Ric(\omega) - i\partial\bar{\partial} \log g(\mu_{\omega})  = \omega \]
where \(\Ric(\omega)\) denotes the Ricci form of \(\omega\), and \(\mu_{\omega}\) denotes the canonically normalized moment map of \(\omega\) with respect to the compact torus action. 
Early results on such metrics, in this generality, were obtained by Mabuchi \cite{Mabuchi_2003}. 
In \cite{Han-Li_2023}, Han and Li proved a version of the Yau-Tian-Donaldson conjecture for weighted solitons in a weak sense, on singular Fano varieties, in terms of the corresponding notion of weighted K-polystability, defined rather in terms of the action of the complexified torus \(T=(\bbC^*)^k\).  
We do not recall the general definition here, but note that it boils down to checking that, for any equivariant, special test configuration, the weighted Donaldson-Futaki invariant is non-negative, and is zero if and only if the test configuration is a product. 
If the latter equality case is not included, the variety is called weighted K-semistable. 
The Donaldson-Futaki invariants associated to product test configurations may further be interpreted as a weighted version of the classical Futaki invariant \cite{Futaki_1983}. 

A Fano toric variety, that is, a Fano variety equipped with an effective regular action of a torus of the same dimension, satisfies a particularly nice property with respect to weighted K-stability. 
Namely, it is weighted K-polystable if and only if it is weighted K-semistable, if and only if its weighted Futaki invariant vanishes. 
This property follows essentially from the fact that the only equivariant special test configurations for toric varieties are product test configurations \cite{Donaldson_2002}, and \cite{Han-Li_2023}. 
However the first results of this nature date back to the characterization of existence of Kähler-Ricci solitons on Fano toric manifolds by Wang and Zhu \cite{Wang-Zhu_2004}. 

In \cite{Delcroix_2020}, we proved that a similar property is satisfied by horospherical \(G\)-varieties. 
A horospherical \(G\)-variety is a variety equipped with an action of a connected reductive group \(G\) such that a \(G\)-stable open set \(X_0\subset X\) is a homogeneous fibre bundle over a generalized flag variety \(G/P\), with fiber a toric variety. 
We proved that for such a variety, the only \(G\)-equivariant special test configurations are product test configurations. 
As a consequence, if \(X\) is a Fano horospherical variety equipped with an action of a torus \(T\) that commutes with the action of \(G\), and \(g\) is a weight function for this \(T\)-action, then \(X\) is \(g\)-K-polystable if and only if its \(g\)-Futaki invariant vanishes. 

Linsheng Wang recently showed that Fano threefolds in families number 2-28 (blowup of \(\bbP^3\) along a plane cubic curve) and 3-14 (further blowup along a point not on the exceptional divisor) exhibit a similar property \cite{Wang_Kpoly,Miao-Wang_moduli}. 
Let \(X\) be a Fano manifold in one of these two families. 
Let \(\bbC^*\) act on \(X\) as the maximal connected reductive subgroup of its automorphism group (in the case of family 3-14, this is actually the full connected automorphism group), and let \(\Delta\) be the canonically normalized moment image of \(c_1(X)\).  
Then for any weight function \(g:\Delta \to \bbR_{>0}\), \(X\) is \(g\)-K-polystable if and only if its \(g\)-Futaki invariant vanishes. 
These examples are of particular interest because in this case, there exists non-product equivariant special test configurations, but they turn out to play no role in weighted K-polystability. 

The purpose of this paper is to investigate other instances of such a behavior via a case study. 
We will mainly study the case of spherical Fano threefolds. 
Apart from the toric examples, the spherical Fano threefolds are (in Mori-Mukai numbering \cite{Mori-Mukai_1981} together with a geometric  description): 
\begin{description}
\item[1-16] the quadric \(Q^3\),
\item[2-29] the blowup of the quadric \(Q^3\) along a conic,
\item[2-30] the blowup of the quadric \(Q^3\) at a point,
\item[2-31] the blowup of \(Q^3\) along a line, 
\item[2-32] the variety \(W\), a divisor in \(\bbP^2\times \bbP^2\) of bidegree \((1,1)\), 
\item[3-19] the blowup of \(Q^3\) along two points, not contained in a line in \(Q^3\), 
\item[3-20] the blowup of \(Q^3\) in the disjoint union of two lines, 
\item[3-22] the blowup of \(\bbP^1\times \bbP^2\) along a conic in a \(\bbP^2\) fiber, 
\item[3-24] a complete intersection of degree \((1,1,0)\) and \((0,1,1)\) in \(\bbP^1\times \bbP^2\times \bbP^2\)
\item[4-7] the blowup of \(W\) along the disjoint union of a curve of degree \((1,0)\) and a curve of degree \((0,1)\), 
\item[4-8] the blowup of \(\bbP^1\times \bbP^1\times \bbP^1\) along a curve of degree \((0,1,1)\). 
\end{description}

To state our main results, we introduce the terminology of \emph{weight-insensitive K-polystability}: a Fano variety \(X\) (or more generally, a \(\bbQ\)-Fano log pair \((X,D)\))  equipped with a regular action of a torus \(T=(\bbC^*)^k\) with canonically normalized moment image \(\Delta\) is weight-insensitive K-polystable if for any weight function \(g:\Delta \to \bbR_{>0}\), \(X\) is \(g\)-K-polystable if and only if its \(g\)-Futaki invariant vanishes. 
We discover several new weight-insensitive K-polystable \(T\)-manifolds.

\begin{thm}
\label{thm:weight-insensitive}
All spherical Fano threefolds but the Fano threefold 2-29 are weight-insensistive K-polystable with respect to the action of the connected center of a Levi subgroup of their automorphism group. 
\end{thm}

Let us stress here that we include the data of the torus action in the definition. 
For the manifolds above, we do not prove that they are weight-insensitive K-polystable with respect to the action of a maximal torus in their automorphism group. 
The latter notion would likely be the most general to hope for, but in this case we cannot use the spherical symmetry since the spherical action does not commute with the maximal torus action. 
Furthermore, for several natural choices of weighted solitons, it is natural to consider only the action of the center of a Levi subgroup of the automorphism group. 
This is the case for Kähler-Ricci solitons, for Mabuchi solitons, as well as for the notion of weighted solitons encoding the existence of Sasaki-Einstein metrics on the link of the cone over the Fano variety. 
Another extreme case would be to consider no torus action. In this case, weight-insensitive K-polystability is equivalent to K-polystability. 

A quick note on Mabuchi solitons: it was proved in \cite{Delcroix-Hultgren_2019} that Fano threefold 2-30 does not admit Mabuchi solitons, and it was proved previously that several toric Fano threefolds do not admit Mabuchi solitons \cite{Yao_2022}. 
This is not in contradiction with the property of weight-insensitive K-polystability, since in these cases the obstruction is that the weight corresponding to Mabuchi solitons is \emph{not} positive. 

To disprove weight-insensitive K-polystability with respect to a maximal torus action, it is enough to disprove it for a weight associated to a lower-dimensional torus action (still with vanishing weighted Futaki invariant). 
Linsheng Wang provided in \cite{Wang_Kpoly} a singular example which is not weight-insensitive K-polystable, the variety obtained by optimal degeneration of some Fano threefolds in family 2-23 \cite{Miao-Wang_optimal}. 
We will provide here simple smooth counterexamples by considering appropriate \(\bbC^*\)-actions on these. 

In \cite{Wang_optimal}, Linsheng Wang highlights an interesting question that further extends this notion of weight-insensitivity. 
Namely, for a K-unstable Fano variety \(X\), he considers variants of the H-functional involved in unweighted optimal degenerations, and constructs associated optimal degenerations which are weighted K-polystable for a weight associated with the functional. 
He asks the question: is this optimal degeneration independent of the choice of variant of the \(H\)-functional (under some assumptions). 
More generally, one can expect the existence of weighted optimal degenerations for arbitrary weight functions. 
This notion would associate to a weighted K-unstable variety for a certain weight \(g\), a certain degeneration which is weighted K-polystable for another weight \(\tilde{g}\), related with \(g\). 
This notion is a generalization of the optimal degeneration from a (unweighted) K-unstable variety to a variety with a (weak) Kähler-Ricci soliton, first studied in the differential geometric setting of the Hamilton-Tian conjecture \cite{Chen-Sun-Wang_2018,Dervan-Szekelyhidi_2020}, then in the algebro-geometric setting \cite{BLXZ_2023}.  
This should be, in general, a two-step process. 
From a \(g\)-K-unstable variety \(X\), one first degenerates to a \(\tilde{g}\)-K-semistable variety \(X_0\), then \(X_0\) degenerates to a \(\tilde{g}\)-K-polystable variety \(Y\). 
The second step is already carried out in full generality in \cite{Han-Li_2024, BLXZ_2023}, where they show the existence and uniqueness of such a \(Y\), starting from a weighted K-semistable \(X_0\). 
Wang's work in \cite{Wang_optimal} constitutes an initial advance in the first step of the process. 
We now state an extension of Wang's question: for which \(X\) equipped with a torus action do we have that the optimal degeneration \(Y\) is independent of the chosen weight. 

Wang's question was originally presented as a conjecture in the first arXiv version of \cite{Wang_optimal}, and was one motivation for the present article. 
Little is known so far on the otimal degenerations of the singular example from \cite{Wang_Kpoly} which is not weight-insensitive K-polystable, so it was not presented as a counterexample to this conjecture. 
For our new examples, we can further characterize the weighted optimal degenerations for certain weights. 

\begin{thm}
\label{thm:weight-sensitive}
The Kähler-Einstein Fano threefold 2-29 is not weight-insensitive K-polystable with respect to the action of the connected center of its automorphism group. 
Furthermore, the quadric threefold \(Q^3\) admits at least one non-trivial \(\bbC^*\)-action for which it is not weight-insensitive K-polystable. 
In both cases, we exhibit an explicit weight \(g\) such that the manifold \(X\) is strictly \(g\)-K-semistable, and its optimal degeneration is a Gorenstein toric variety not isomorphic to \(X\). 
\end{thm}

The basic idea behind these counterexamples is that, by \cite{Delcroix_2020,Delcroix_EpiGA}, the spherical threefolds considered admit a unique \(G\)-equivariant special test configuration up to twist, where \(G=\SL_2(\bbC)\times \bbC^*\) in our case. 
Furthermore, by the criterion from \cite{Delcroix_2020} as generalized to weighted K-polystability in \cite{Li-Li-Wang}, the central fiber of this test configuration is weighted K-semistable, and horospherical, hence weighted K-polystable. 
As a consequence, the optimal degeneration can only be induced by that test configuration. 

It is actually easy to describe the special test configuration for our examples. 
Consider the quadric \(Q^3\) as defined by the equation 
\[ x_0x_2-x_1^2 + x_3x_4 = 0 \]
in \(\bbP^4\) with homogeneous coordinates \([x0:\cdots:x_4]\). 
Then the test configuration \(\mathcal{X}\) defined by the equation 
\[ x_0x_2-x_1^2 + zx_3x_4 = 0 \]
in \(\bbC\times \bbP^4\) (with coordinate \(z\) for the \(\bbC\) factor) is a \(G\)-equivariant special non-product test configuration, with central fiber the horospherical \(G\)-variety which is the singular quadric defined by the equation 
\[ x_0x_2-x_1^2 = 0 \]
in \(\bbP^4\). 
This singular quadric is further toric and Gorenstein, corresponding to the Fano polytope with reflexive ID 2 in the Graded Ring Database \cite{grdb}. 
The associated moment polytope is the simplex obtained as the convex hull of the four points \((-1,-1,-1)\), \((5,-1,-1)\), \((-1,2,-1)\) et \((-1,-1,2)\) in \(\bbR^3\). 

For the threefold 2-29, we construct the special test configuration by blowing up the subvariety 
\(\mathcal{Z}\subset \mathcal{X}\) defined by \(x_3=x_4=0\) in the test configuration \(\mathcal{X}\) for \(Q^3\). 
Again, the central fiber is \(G\)-horospherical but as well toric and Gorenstein, with reflexive ID 19 in \cite{grdb}. 

Let us note that our examples are not (unweighted) K-unstable, and the weights considered in our examples do not quite satisfy the added hypotheses in \cite{Wang_optimal} (that some primitive of the weight is log-convex). 
This does not seem like a major issue to us, since we produce strictly weighted K-semistable examples, for which the assumption on the weights is not necessary to define an optimal degeneration in \cite{Wang_optimal}. 
Strictly speaking though, we do not provide a counterexample to the original conjecture in \cite{Wang_optimal}. 

If one switches to the log Fano setting, we can provide an (unweighted) K-unstable example. We illustrate this by considering log Fano pairs with underlying manifold the threefold 2-29. Let \(X\) denote the manifold 2-29, and let \(E\) denote the exceptional divisor of the blowup map \(X\to Q^3\). 

\begin{thm}
\label{thm:logpairs}
For \(t_0=\frac{\sqrt{10}-2}{3}\), the pair \((X,t_0E)\) is strictly (unweighted) K-semistable. 
Furthermore, there are smooth log concave weights for which \((X,t_0E)\) is weighted K-polystable.
\end{thm}

Finally, we observe that the above results carry through to higher dimensional quadrics. 
Under the action of \(\SO_{n-2}(\bbC)\times \SO_2(\bbC)\subset \SO_{n+2}(\bbC)\) on the quadric \(Q^{n-2}\), 
there is still a unique equivariant special test configuration up to twist. 
The candidate optimal degeneration is in fact still the degeneration to the singular quadric 
\[ x_0^2+\cdots +x_{n-3}^2 = 0 \]
in \(\bbP^{n-1}\). 
This singular quadric is no longer toric in general, but it is still rank two horospherical. 

\begin{thm}
\label{thm:high-dim}
For any \(n\geq 5\), the quadric \(Q^{n-2}\) is not weight-insensitive K-polystable. 
There exist a weight for which it is strictly weighted K-semistable, and its weighted optimal degeneration is the singular quadric with equation 
\(x_0^2+\cdots +x_{n-3}^2=0\) 
in \(\bbP^{n-1}\). 
\end{thm}

More generally, a spherical Fano variety admits a unique  special test configuration up to twist if and only if its valuation cone is a half space. These constitute in general natural candidates to study with respect to weighted optimal degenerations. 

\begin{ack}
The author is partially funded by ANR-21-CE40-0011 JCJC project MARGE.
We thank Linsheng Wang and Tiago Duarte Guerreiro for conversations on \cite{Wang_optimal} and on the present article. 
\end{ack}

\section{Recollections on spherical varieties}

\subsection{Combinatorial data associated with spherical varieties}

We first recall some of the basics of the theory of spherical varieties that will be needed. 
We refer to \cite{Knop_1991,Brion_1989,Delcroix-Montagard} for more details. 
Let \(X\) be a normal variety, and let \(G\) be a connected reductive group, acting regularly on \(X\). 
If a Borel subgroup of \(G\) acts with an open orbit in \(X\), then \(X\) is called a spherical \(G\)-variety. 
We assume for the rest of this section that \(X\) is a spherical \(G\)-variety. 
We further choose a base point \(x\) in the open orbit of \(G\), we denote by \(H\) its stabilizer in \(G\), and we choose a Borel subgroup \(B\) such that \(BH\) is open in \(G\) (in other words, such that \(B\cdot x\) is open in \(X\)). 

The \emph{weight lattice} of \(X\), denoted by \(M\), is the set of all characters \(m\) of \(B\) such that there exists a rational function \(f_m\in \bbC(X)\) with \(b\cdot f_m = m(b) f\).  
Since \(B\) acts with an open orbit, such an \(f\) is uniquely determined by \(m\) up to a multiplicative constant. 
We denote by \(X^*(B)\) the set of characters of \(B\), so that \(M\subset X^*(B)\). 
The set \(M\) is a torsion-free finite rank abelian group isomorphic to some \(\bbZ^r\). 
The integer \(r\) is called the rank of the \(G\)-variety \(X\). 
We denote by \(N=\Hom(M,\bbZ)\) the dual abelian group. 

Any \(\bbR\)-valued valuation of \(\bbC(X)\) induces an element of \(N\otimes \bbR=\Hom(M,\bbR)\) by restriction to \(B\)-eigenfunctions as above. 
The set \(\mathcal{V}\) of all the elements of \(N\otimes \bbR\) obtained by restriction of a \(G\)-invariant valuations forms a convex polyhedral cone in \(N\otimes \bbR\) called the \emph{valuation cone}. 

Assume that \(X\) is \(\bbQ\)-Fano. Then there is a natural linearization of \(K_X^{-m}\) for \(m\) divisible enough, and from this a natural moment polytope \(\Delta^+(X)\subset X^*(B)\otimes \bbR\). 
This is by definition the closure of the set of all \(\frac{\lambda}{m}\) where \(K_X^{-m}\) is a line bundle, \(\lambda\in X^*(B)\), and there exist a section \(s\in H^0(X,K_X^{-m})\) such that for all \(b\in B\), \(b\cdot s = \lambda(b)s\) for the action of \(G\) induced by the natural linearization of \(K_X^{-m}\). 

Choose a maximal torus \(T\subset B\), and let \(R\) and \(R^+\) denote the set of roots and the set of positive roots of \(G\) with respect to these choices. 
Let \(R_X^+\subset R^+\) denote the subset of positive roots \(\alpha\) such that \(\left\{\alpha,\cdot\right\}\) does not vanish identically on \(\Delta^+(X)\), where \(\left\{\cdot,\cdot\right\}\) denotes the Killing form. 
Then we call Duistermaat-Heckman polynomial of \(X\) the polynomial \(P_{DH}\) on \(X^*(B)\otimes \bbR\) defined by 
\[ P_{DH}(p) = \prod_{\alpha\in R_X^+}\frac{\left\{\alpha, p\right\}}{\left\{\alpha,\varpi\right\}}\]
where \(\varpi=\frac{1}{2}\sum_{\alpha\in R^+}\alpha\). 

\subsection{Weighted K-stability of spherical Fano varieties}

We recall first the classification of equivariant special test configurations for spherical varieties from \cite{Delcroix_2020}. 
The additional characterization of special test configurations up to twist follows from \cite{Delcroix_EpiGA}. 

\begin{thm}[\cite{Delcroix_2020,Delcroix_EpiGA}]
\label{thm_degenerations}
Let \(X\) be a \(\bbQ\)-Fano variety, spherical under the action of a connected reductive group \(G\). 
Then \(G\)-equivariant special test configurations for \(X\) are in bijection with rational rays in the valuation cone \(\mathcal{V}(X)\). 
Furthermore, \(G\)-equivariant special test configurations up to twists are in bijection with rational rays in the projection of the valuation to the quotient \(\frac{N\otimes \bbR}{\mathcal{V}\cap -\mathcal{V}}\). 
\end{thm}

Li, Li and Wang proved in \cite[Theorem~1.3]{Li-Li-Wang} a criterion for weighted K-polystability of spherical Fano varieties, a generalization of the main result of \cite{Delcroix_2020} which dealt only explicitly with the case of Kähler-Ricci soliton weights. 
We reformulate it as follows. 

First, let us stress that we focus on torus actions that commute with the action of \(G\). 
It is known (\cite{Losev_2009}, see the discussion in \cite[section~3.1.3]{Delcroix_2020}) that, for a spherical \(G\)-variety \(X\), the connected component \(\bbT\) of the group \(\Aut^G(X)\) of \(G\)-equivariant automorphisms of \(X\) is a torus. 
Let us assume, to simplify notations, that the action of \(G\) is almost-faithful (that is, faithful up to a finite central subgroup), and that the image of \(G\) in \(\Aut(X)\) contains \(\bbT\). 
Then it is furthermore known that the moment polytope of \(X\) with respect to the action of \(\bbT\) is the image of \(\Delta^+\) under the map 
\[ \pi_{\bbT}: X^*(B)\otimes \bbR \to X^*(T/([G,G]\cap T))\otimes \bbR = X^*(\bbT)\otimes \bbR \] 
the latter identification being obtained via the isogenies between \(\bbT\), \(Z(G)^0\) and \(T/([G,G]\cap T)\). 

\begin{thm}[\cite{Delcroix_2020} and \cite{Li-Li-Wang}]
\label{thm_weighted_K-stab_spherical}
The \(\bbQ\)-Fano spherical variety \(X\) is \(G\)-equivariantly \(g\)-K-semistable if and only if for any \(\xi\in \mathcal{V}\), 
\[ \int_{\Delta^+} \langle p-\kappa, \xi \rangle g(\pi_{\bbT}(p)) P_{DH}(p) \mathop{dp} \leq 0 \]
where \(\kappa = \sum_{\alpha\in R_X^+} \alpha\). 
If furthermore, equality holds if and only if \(\xi\in \mathcal{V}\cap -\mathcal{V}\), then \(X\) is \(G\)-equivariantly \(g\)-K-polystable. 
The vanishing of the above integral for all \(\xi\in \mathcal{V}\cap -\mathcal{V}\) is equivalent to the vanishing of the weighted Futaki invariant. 
\end{thm}

In the above statement, \(\mathop{dp}\) denotes a Lebesgue measure on the affine space \(\kappa + M\otimes \bbR\). The condition does not depend on a particular choice of such a Lebesgue measure, although it is customary to choose the Lebesgue measure normalized by the lattice \(M\), so that, as a sanity check, one can check consistency with the formula for the anticanonical degree: 
\[ (-K_X)^n = (\dim X) ! \int_{\Delta^+} P_{DH}(p) \mathop{dp}  \]
We will not choose this particular normalization in the following. 

\begin{rem}
Although the variant is not written in the literature, we note that the above statement works as well with the same proof for a klt log Fano pair \((X,tD)\) where \(X\) is a \(G\)-spherical Fano variety, \(D\) is a \(G\)-stable divisor, and \(t\in \bbR\). 
It suffices to replace the moment polytope \(\Delta^+\) in the statement by the moment polytope of the anticanonical divisor of the pair \(-K_X-tD\). 
\end{rem}

\section{Spherical Fano threefolds}

\subsection{Faithful spherical actions on Fano threefolds}

Since toric Fano varieties are weight-insensitive K-polystable with respect to any torus action, we focus on non-toric spherical Fano threefolds. 
In Pasquier's \cite{Pasquier_2006} and Hofscheier's combined PhD theses (the latter is unavailable online), the classification of all possible faithful spherical actions of a connected reductive group on a Fano threefolds was essentially achieved. These results were reproved and generalized in \cite{Delcroix-Montagard}. 
From these we gather: 

\begin{thm}{\cite{Pasquier_2006,Delcroix-Montagard}}
In Mori-Mukai numbering \cite{Mori-Mukai_1981}, Fano threefolds number 1-16, 2-29, 2-30, 2-31, 2-32, 3-19, 3-20, 3-22, 3-24, 4-7 and 4-8 are the only non-toric, spherical Fano threefolds. 
All but 1-16 and 2-32 admit only one faithful spherical action up to equivariant isomorphism, induced by an almost-faithful action of \(\SL_2\times \bbG_m\). 
For 1-16 there are two different faithful spherical actions of rank two induced by \(\SL_2\times \bbG_m\)-actions, there are also rank one and rank zero spherical actions that we will not consider here. 
For 2-32 the only other faithful spherical action is the rank zero action of its automorphism group. 
We summarize the correspondence with the identifiers from \cite{Delcroix-Montagard} in Table~\ref{table:spherical3folds}, for rank two faithful spherical actions.  
\end{thm}

\begin{table}
    \centering
    \begin{tabular}{ccccccc}
      Mori-Mukai& 1-16 & 1-16 & 2-29 & 2-30 & 2-31 & 2-32 \\
       \midrule
        Delcroix-Montagard & 3-2-4 & 3-2-18 & 3-2-19 & 3-2-21 & 3-2-6 & 3-2-5 \\
        \\
       Mori-Mukai & 3-19 & 3-20 & 3-22 & 3-24 & 4-7 & 4-8 \\
       \midrule
        Delcroix-Montagard & 3-2-23 & 3-2-9 & 3-2-17 & 3-2-8 & 3-2-11 & 3-2-3 \\
        \\
    \end{tabular}
    \caption{Faithful \(\SL_2\times \bbG_m\)-spherical actions on Fano threefolds}
    \label{table:spherical3folds}
\end{table}

\subsection{Weighted K-stability}

The goal of this subsection is to prove Theorem~\ref{thm:weight-insensitive}. 
The result is known for toric manifolds. 
It is also true for the quadric \(Q^3\): its automorphism group is semisimple, so its connected center is trivial. Since \(Q^3\) is K-polystable, it is weight-insensitive K-polystable with respect to the trivial torus action. We will see later that it is not always true for non-trivial \(\bbC^*\)-actions. 

In view of Theorem~\ref{thm_weighted_K-stab_spherical}, it suffices to know the valuation cone and the moment polytope \(\Delta^+\subset X^*(B)\otimes \bbR\) of a spherical variety to determine when it is weighted K-polystable. 
The combinatorial data given in \cite{Delcroix-Montagard} is enough to determine these as explained in \cite[Section~2.8]{Delcroix-Montagard}. 
In order to identify the optimal degeneration, we will further indicate the lattice \(M\subset X^*(B)\) for each example (or rather, we will indicate the points in the intersection of the translated lattice \(\kappa+M\) with \(\Delta^+\)). 

\begin{figure}
\centering
\begin{tikzpicture}[scale=0.6]
\draw[dotted] (0,3) grid (6,-3);
\draw (0,0) node[left]{0};
\draw (0,1) node[left]{\(\chi\)};
\draw[thick] (0,-3) -- (6,0) -- (0,3) -- cycle;
\draw (2,0) node[below left]{\(\alpha\)};
\draw (0,3) node{\(\bullet\)}; \draw (0,1) node{\(\bullet\)}; \draw (0,-1) node{\(\bullet\)}; \draw (0,-3) node{\(\bullet\)}; \draw (2,0) node{\(\bullet\)}; \draw (2,-2) node{\(\bullet\)}; \draw (2,2) node{\(\bullet\)}; \draw (4,1) node{\(\bullet\)}; \draw (4,-1) node{\(\bullet\)}; \draw (6,0) node{\(\bullet\)};
\draw (3,-3) node[below]{3-2-18};
\end{tikzpicture}
\begin{tikzpicture}[scale=0.6]
\draw[dotted] (0,3) grid (5,-3);
\draw (0,0) node[left]{0};
\draw (0,1) node[left]{\(\chi\)};
\draw[thick] (0,3) -- (4,1) -- (4,-1) -- (0,-3) -- cycle;
\draw (2,0) node[below left]{\(\alpha\)};
\draw (0,3) node{\(\bullet\)}; \draw (0,1) node{\(\bullet\)}; \draw (0,-1) node{\(\bullet\)}; \draw (0,-3) node{\(\bullet\)}; \draw (2,0) node{\(\bullet\)}; \draw (2,-2) node{\(\bullet\)}; \draw (2,2) node{\(\bullet\)}; \draw (4,1) node{\(\bullet\)}; \draw (4,-1) node{\(\bullet\)}; 
\draw (3,-3) node[below]{3-2-19};
\end{tikzpicture}
\begin{tikzpicture}[scale=0.6]
\draw[dotted] (0,3) grid (6,-3);
\draw (0,0) node[left]{0};
\draw (0,1) node[left]{\(\chi\)};
\draw[thick] (0,3) -- (6,0) -- (4,-1) -- (0,-1) -- cycle;
\draw (2,0) node[below left]{\(\alpha\)};
\draw (0,3) node{\(\bullet\)}; \draw (0,1) node{\(\bullet\)}; \draw (0,-1) node{\(\bullet\)};  \draw (2,0) node{\(\bullet\)}; \draw (2,2) node{\(\bullet\)}; \draw (4,1) node{\(\bullet\)}; \draw (4,-1) node{\(\bullet\)}; \draw (6,0) node{\(\bullet\)};
\draw (3,-3) node[below]{3-2-21};
\end{tikzpicture}
\begin{tikzpicture}[scale=0.6]
\draw[dotted] (0,3) grid (6,-3);
\draw (0,0) node[left]{0};
\draw (0,1) node[left]{\(\chi\)};
\draw[thick] (0,-1) -- (4,-1) -- (6,0) -- (4,1) -- (0,1) -- cycle;
\draw (2,0) node[below left]{\(\alpha\)};
\draw (0,1) node{\(\bullet\)}; \draw (0,-1) node{\(\bullet\)}; \draw (2,0) node{\(\bullet\)}; \draw (4,1) node{\(\bullet\)}; \draw (4,-1) node{\(\bullet\)}; \draw (6,0) node{\(\bullet\)};
\draw (3,-3) node[below]{3-2-23};
\end{tikzpicture}
\begin{tikzpicture}[scale=0.6]
\draw[dotted] (0,3) grid (6,-3);
\draw (0,0) node[left]{0};
\draw (0,1) node[left]{\(\chi\)};
\draw[thick] (0,-1) -- (2,-1) -- (6,-1) -- (6,0) -- (2,1) -- (0,1) -- cycle;
\draw (2,0) node[below left]{\(\alpha\)};
\draw (2,0) node{\(\bullet\)}; 
\draw (2,1) node{\(\bullet\)}; \draw (2,-1) node{\(\bullet\)}; \draw (6,0) node{\(\bullet\)}; \draw (6,-1) node{\(\bullet\)}; 
\draw (3,-3) node[below]{3-2-17};
\end{tikzpicture}
\begin{tikzpicture}[scale=0.6]
\draw[dotted] (0,3) grid (5,-3);
\draw (0,0) node[left]{0};
\draw (0,1) node[left]{\(\chi\)};
\draw[thick] (0,-1) -- (4,-1) -- (4,0) -- (2,1) -- (0,1) -- cycle;
\draw (2,0) node[below left]{\(\alpha\)};
\draw (2,0) node{\(\bullet\)}; 
\draw (2,1) node{\(\bullet\)}; \draw (2,-1) node{\(\bullet\)}; \draw (4,0) node{\(\bullet\)}; \draw (4,-1) node{\(\bullet\)}; \draw (0,1) node{\(\bullet\)}; \draw (0,-1) node{\(\bullet\)}; \draw (0,0) node{\(\bullet\)};
\draw (3,-3) node[below]{3-2-3};
\end{tikzpicture}
\begin{tikzpicture}[scale=0.6]
\draw[dotted] (0,3) grid (4,-3);
\draw (0,0) node[left]{0};
\draw (0,1) node[left]{\(\chi\)};
\draw[thick] (0,0) -- (3,-3) -- (3,3) -- cycle;
\draw (2,0) node[below left]{\(\alpha\)};
\draw (2,0) node{\(\bullet\)}; 
\draw (0,0) node{\(\bullet\)}; \draw (1,-1) node{\(\bullet\)}; \draw (1,1) node{\(\bullet\)}; \draw (2,-2) node{\(\bullet\)}; \draw (2,2) node{\(\bullet\)}; \draw (3,-1) node{\(\bullet\)}; \draw (3,1) node{\(\bullet\)}; \draw (3,-3) node{\(\bullet\)}; \draw (3,3) node{\(\bullet\)};
\draw (2,-3) node[below]{3-2-4};
\end{tikzpicture}
\begin{tikzpicture}[scale=0.6]
\draw[dotted] (0,3) grid (4,-3);
\draw (0,0) node[left]{0};
\draw (0,1) node[left]{\(\chi\)};
\draw[thick] (0,0) -- (2,-2) -- (4,0) -- (2,2) -- cycle;
\draw (2,0) node[below left]{\(\alpha\)};
\draw (2,0) node{\(\bullet\)}; 
\draw (0,0) node{\(\bullet\)}; \draw (1,-1) node{\(\bullet\)}; \draw (1,1) node{\(\bullet\)}; \draw (2,-2) node{\(\bullet\)}; \draw (2,2) node{\(\bullet\)}; \draw (3,-1) node{\(\bullet\)}; \draw (3,1) node{\(\bullet\)}; \draw (4,0) node{\(\bullet\)}; 
\draw (2,-3) node[below]{3-2-5};
\end{tikzpicture}
\begin{tikzpicture}[scale=0.6]
\draw[dotted] (0,3) grid (4,-3);
\draw (0,0) node[left]{0};
\draw (0,1) node[left]{\(\chi\)};
\draw[thick] (0,0) -- (2,-2) -- (3,-1) -- (3,3) -- cycle;
\draw (2,0) node[below left]{\(\alpha\)};
\draw (2,0) node{\(\bullet\)}; 
\draw (0,0) node{\(\bullet\)}; \draw (1,-1) node{\(\bullet\)}; \draw (1,1) node{\(\bullet\)}; \draw (2,-2) node{\(\bullet\)}; \draw (2,2) node{\(\bullet\)}; \draw (3,-1) node{\(\bullet\)}; \draw (3,1) node{\(\bullet\)}; \draw (3,3) node{\(\bullet\)};
\draw (2,-3) node[below]{3-2-6};
\end{tikzpicture}
\begin{tikzpicture}[scale=0.6]
\draw[dotted] (0,3) grid (4,-3);
\draw (0,0) node[left]{0};
\draw (0,1) node[left]{\(\chi\)};
\draw[thick] (0,0) -- (1,-1) -- (3,-1) -- (4,0) -- (2,2) -- cycle;
\draw (2,0) node[below left]{\(\alpha\)};
\draw (2,0) node{\(\bullet\)}; 
\draw (0,0) node{\(\bullet\)}; \draw (1,-1) node{\(\bullet\)}; \draw (1,1) node{\(\bullet\)};  \draw (2,2) node{\(\bullet\)}; \draw (3,-1) node{\(\bullet\)}; \draw (3,1) node{\(\bullet\)}; \draw (4,0) node{\(\bullet\)}; 
\draw (2,-3) node[below]{3-2-8};
\end{tikzpicture}
\begin{tikzpicture}[scale=0.6]
\draw[dotted] (0,3) grid (4,-3);
\draw (0,0) node[left]{0};
\draw (0,1) node[left]{\(\chi\)};
\draw[thick] (0,0) -- (2,-2) -- (3,-1) -- (3,1) -- (2,2) -- cycle;
\draw (2,0) node[below left]{\(\alpha\)};
\draw (2,0) node{\(\bullet\)}; 
\draw (0,0) node{\(\bullet\)}; \draw (1,-1) node{\(\bullet\)}; \draw (1,1) node{\(\bullet\)}; \draw (2,-2) node{\(\bullet\)}; \draw (2,2) node{\(\bullet\)}; \draw (3,-1) node{\(\bullet\)}; \draw (3,1) node{\(\bullet\)}; 
\draw (2,-3) node[below]{3-2-9};
\end{tikzpicture}
\begin{tikzpicture}[scale=0.6]
\draw[dotted] (0,3) grid (4,-3);
\draw (0,0) node[left]{0};
\draw (0,1) node[left]{\(\chi\)};
\draw[thick] (0,0) -- (1,-1) -- (3,-1) -- (4,0) -- (3,1) -- (1,1) -- cycle;
\draw (2,0) node[below left]{\(\alpha\)};
\draw (2,0) node{\(\bullet\)}; 
\draw (0,0) node{\(\bullet\)}; \draw (1,-1) node{\(\bullet\)}; \draw (1,1) node{\(\bullet\)};  \draw (3,-1) node{\(\bullet\)}; \draw (3,1) node{\(\bullet\)}; \draw (4,0) node{\(\bullet\)}; 
\draw (2,-3) node[below]{3-2-11};
\end{tikzpicture}
\caption{Moment polytopes of non-toric \(\SL_2\times\bbC^*\)-spherical Fano threefolds}
\label{fig:polytope_DH}
\end{figure}
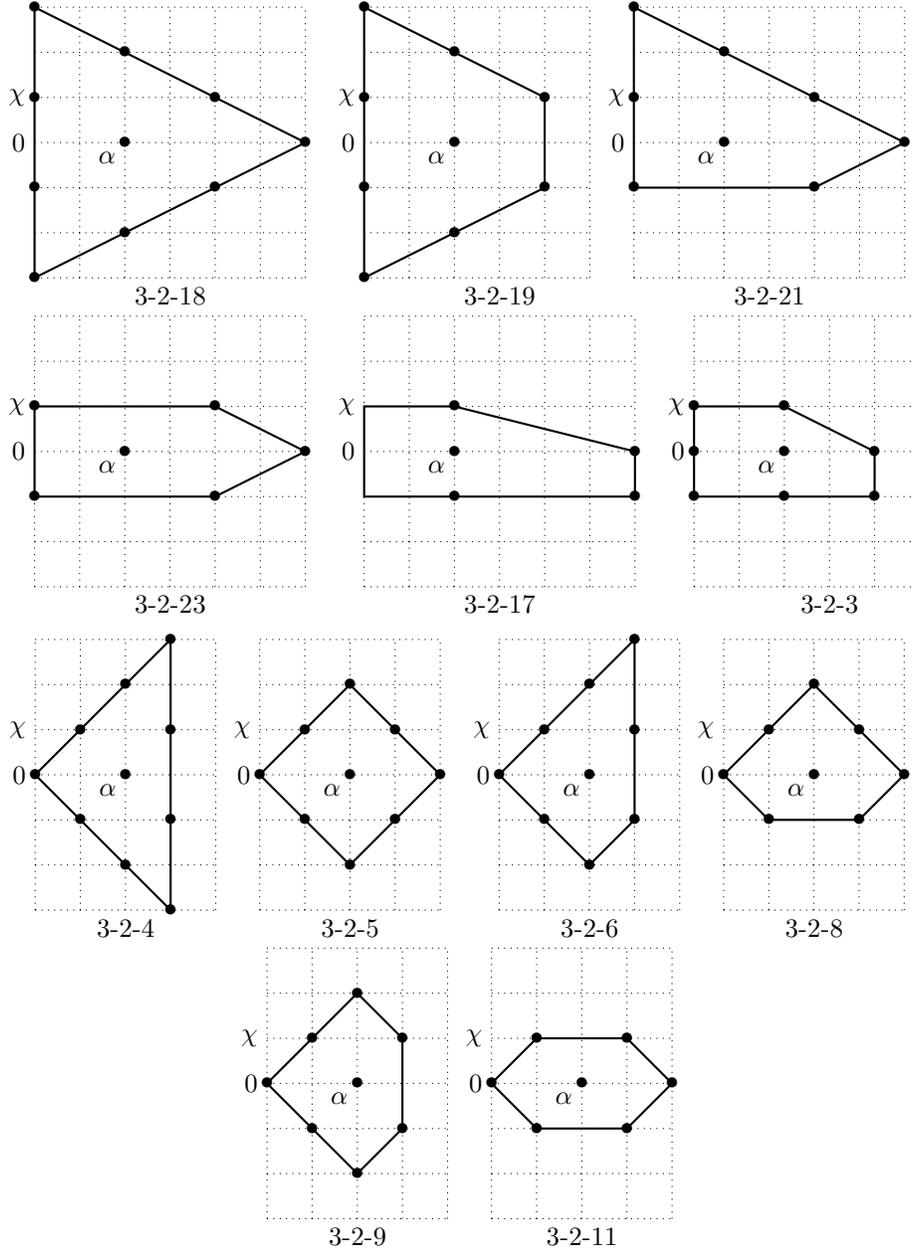

Here we are working with \(G=\SL_2\times \bbG_m\). 
Fix a maximal torus \(T\subset B\), and let \(\alpha\in X^*(B)\) be the unique positive root. 
Let \(\chi\) be a primitive generator of \(X^*(T/([G,G]\cap T))\). 
The moment polytopes for all the (almost-)faithful spherical actions of \(G\) on non-toric Fano threefolds are depicted in Figure~\ref{fig:polytope_DH}. 

The spherical root in each case is the positive primitive multiple of \(\alpha\) in \(M\), so that the valuation cone is \(\mathcal{V}=\{ \nu\in N\otimes \bbR \mid \nu(\alpha)\leq 0 \}\), 
and \(\mathcal{V}\cap -\mathcal{V} = \bbR \chi^*\) in the dual basis \((\alpha^*,\chi^*)\). 
The Duistermaat-Heckman will not depend on the considered threefolds because in all cases we will have \(R^+=R^+_X=\{\alpha\}\). 
We thus have 
\[ P_{DH}(x\frac{\alpha}{2}+y\chi) = \frac{x}{4} \]
Furthermore, in all cases, we have 
\[ \kappa = \alpha \]

We identify \(X^*(B)\otimes \bbR\) with \(\bbR^2\) by using the coordinates \((x,y)\) in the basis \((\frac{\alpha}{2},\chi)\). 
A weight \(g\) is then identified with a function of \(y\), and the criterion for weighted K-polystability reads (up to the choice of the Lebesgue measure) as the two conditions 
\begin{equation*} 
\label{gFut_DH}
\int_{\Delta^+} yg(y)x\mathop{dxdy} = 0 
\end{equation*}
and
\begin{equation*}
\label{stab_DH}
\int_{\Delta^+} (x-2)g(y)x\mathop{dxdy} > 0 
\end{equation*}
We may consider the above expressions as signed measures evaluated at \(g\) by setting  
\[ \nu(g) = \int_{\Delta^+} yg(y)x\mathop{dxdy} \qquad \qquad \mu(g) = \int_{\Delta^+} (x-2)g(y)x\mathop{dxdy} \] 
We compute the measures \(\mu\) and \(\nu\) for all cases and gather the results in Table~\ref{tab:mu_and_nu}.

\begin{table}
\centering
\begin{tabular}{ll}
\toprule
Threefold & Signed measures \\
\midrule
3-2-3 & \(\frac{\mu}{dy}=\frac{16}{3}\indicator{[-1,0]}+\frac{4}{3}(1-2y)(2-y)^2\indicator{[0,1]}\)\\
 & \(\frac{\nu}{dy}=8y\indicator{[-1,0]}+2y(2-y)^2\indicator{[0,1]}\)\\
\midrule 
3-2-4 & \(\frac{\mu}{dy}=\frac{1}{3}y^2(y+3)\indicator{[-3,0]} + \frac{-1}{3}y^2(y-3)\indicator{[0,3]}\)\\
 & \(\frac{\nu}{dy}=\frac{1}{2}y(3-y)(y+3)\indicator{[-3,3]}\)\\
\midrule
3-2-5 & \(\frac{\mu}{dy} = \frac{2}{3}(y+2)^3\indicator{[-2,0]} + \frac{-2}{3}(y-2)^3\indicator{[0,2]}\) \\ 
& \(\frac{\nu}{dy} = 4y(2+y)\indicator{[-2,0]} + 4y(2-y)\indicator{[0,2]}\) \\
\midrule 
3-2-6 & \(\frac{\mu}{dy}=\frac{2}{3}(y+2)^3\indicator{[-2,-1]} +\frac{1}{3}y^2(y+3)\indicator{[-1,0]} + \frac{-1}{3}y^2(y-3)\indicator{[0,3]}\)\\
 & \(\frac{\nu}{dy}=4y(2+y)\indicator{[-2,-1]} +\frac{1}{2}y(3-y)(y+3)\indicator{[-1,3]}\)\\
\midrule 
3-2-8 & \(\frac{\mu}{dy} = \frac{2}{3}(y+2)^3\indicator{[-1,0]} + \frac{-2}{3}(y-2)^3\indicator{[0,2]}\) \\ 
& \(\frac{\nu}{dy} = 4y(2+y)\indicator{[-1,0]} + 4y(2-y)\indicator{[0,2]}\) \\
\midrule
3-2-9 & \(\frac{\mu}{dy} = \frac{2}{3}(y+2)^3\indicator{[-2,-1]} + \frac{1}{3}y^2(y+3)\indicator{[-1,0]} + \frac{-1}{3}y^2(y-3)\indicator{[0,1]} + \frac{-2}{3}(y-2)^3\indicator{[1,2]}\) \\ 
& \(\frac{\nu}{dy} = 4y(2+y)\indicator{[-2,-1]} + \frac{1}{2}y(3-y)(y+3)\indicator{[-1,1]} + 4y(2-y)\indicator{[1,2]}\) \\
\midrule
3-2-11 & \(\frac{\mu}{dy} = \frac{2}{3}(y+2)^3\indicator{[-1,0]} + \frac{-2}{3}(y-2)^3\indicator{[0,1]}\) \\ 
& \(\frac{\nu}{dy} = 4y(2+y)\indicator{[-1,0]} + 4y(2-y)\indicator{[0,1]}\) \\
\midrule 
3-2-17 & \(\frac{\mu}{dy}=36\indicator{[-1,0]}+\frac{4}{3}(3-2y)^2(3-4y)\indicator{[0,1]}\)\\
 & \(\frac{\nu}{dy}=18y\indicator{[-1,0]}+2y(3-2y)^2\indicator{[0,1]}\)\\
\midrule 
3-2-18 & \(\frac{\mu}{dy}=\frac{4}{3}(y+3)^2(2y+3)\indicator{[-3,0]}+\frac{4}{3}(y-3)^2(3-2y)\indicator{[0,3]}\)\\
 & \(\frac{\nu}{dy}=2y(3+y)^2\indicator{[-3,0]}+2y(3-y)^2\indicator{[0,3]}\)\\
\midrule 
3-2-19 & \(\frac{\mu}{dy}=\frac{4}{3}(y+3)^2(2y+3)\indicator{[-3,-1]}+\frac{16}{3}\indicator{[-1,1]}+\frac{4}{3}(y-3)^2(3-2y)\indicator{[1,3]}\)\\
 & \(\frac{\nu}{dy}=2y(3+y)^2\indicator{[-3,-1]}+8y\indicator{[-1,1]}+2y(3-y)^2\indicator{[1,3]}\)\\
\midrule 
3-2-21 & \(\frac{\mu}{dy}=\frac{4}{3}(y+3)^2(2y+3)\indicator{[-1,0]}+\frac{4}{3}(y-3)^2(3-2y)\indicator{[0,3]}\)\\
 & \(\frac{\nu}{dy}=2y(3+y)^2\indicator{[-1,0]}+2y(3-y)^2\indicator{[0,3]}\)\\
\midrule 
3-2-23 & \(\frac{\mu}{dy}=\frac{4}{3}(y+3)^2(2y+3)\indicator{[-1,0]}+\frac{4}{3}(y-3)^2( 3-2y)\indicator{[0,1]}\)\\
 & \(\frac{\nu}{dy}=2y(3+y)^2\indicator{[-1,0]}+2y(3-y)^2\indicator{[0,1]}\)\\
 \midrule
\end{tabular}
\caption{The signed measures \(\mu\) and \(\nu\)}
\label{tab:mu_and_nu}
\end{table}

Let us quickly explain the details of the computation in the example of 3-2-3: we have 
\begin{align*}
\nu(g) & = \int_{\Delta^+} yg(y)x\mathop{dxdy} \\
 & = \int_{y=-1}^0 g(y) y \int_{x=0}^4x\mathop{dx} \mathop{dy}+\int_{y=0}^1 g(y) y \int_{x=0}^{4-2y} x\mathop{dx} \mathop{dy} \\
 & = \int_{y=-1}^0 g(y) 8y \mathop{dy}+\int_{y=0}^1 g(y) 2 y(2-y)^2 \mathop{dy}
\end{align*}
and 
\begin{align*}
\mu(g) & = \int_{\Delta^+} (x-2)g(y)x\mathop{dxdy} \\
 & = \int_{y=-1}^0 g(y) \int_{x=0}^4x(x-2)\mathop{dx} \mathop{dy}+\int_{y=0}^1 g(y) \int_{x=0}^{4-2y} x(x-2)\mathop{dx} \mathop{dy} \\
 & = \int_{y=-1}^0 g(y) \frac{16}{3} \mathop{dy}+\int_{y=0}^1 g(y) \frac{4}{3}(y-2)^2(1-2y) \mathop{dy}
\end{align*}

In these terms, we easily observe that weight-insensitive K-polystability is implied by the (stronger) condition that \(\mu\) is a positive measure with support \(\Delta_{\bbT}\). 
This condition is satisfied for spherical threefolds number 3-2-4, 3-2-5, 3-2-6, 3-2-8, 3-2-9, 3-2-11, and 3-2-23. 

It is not actually necessary for \(\mu\) to be a positive measure. 
For the remaining weight insensitive K-polystable cases, we will rather prove that \(\mu+\lambda\nu\) is a positive measure with support \(\Delta_{\bbT}\) for some well-chosen \(\lambda\in \bbR\). 
This is enough to guarantee weight-insensitive K-polystability since for this notion we only consider weights \(g\) such that \(\nu(g)=0\). 

For 3-2-3, we have the positive measure
\[ \mu+\frac{2}{3}\nu = \left(\frac{16}{3}(1+y)\indicator{[-1,0]} + \frac{4}{3}(2-y)^2(1-y)\indicator{[0,1]}\right) dy \]
For 3-2-17, we have 
\[ \mu+2\nu = \left(36(1+y)\indicator{[-1,0]} + \frac{4}{3}(3-2y)^2(4-y)\indicator{[0,1]}\right)dy \]
And for 3-2-21, 
\[ \mu+\frac{2}{3}\nu = \left(4(y+3)^2(y+1)\indicator{[-1,0]} + \frac{4}{3}(3-y)^3\indicator{[0,3]}\right)dy \]
This finishes the proof of Theorem~\ref{thm:weight-insensitive}. 

Note that we have found a non-trivial \(\bbC^*\)-action on the quadric \(Q^3\) such that \(Q^3\) is weight-insensitive K-polystable with respect to this action. 
In the next section, we prove that it is not true for another choice of \(\bbC^*\)-action. 

\subsection{The quadric threefold and its toric optimal degeneration}

We now focus on the spherical action 3-2-18 on the Fano threefold 1-16, the quadric \(Q^3\). 
In this section, we set \(X=Q^3\), and \(G=\SL_2\times \bbC^*\) acting on \(X\) via the spherical action 3-2-18. 

\subsubsection{Unstable and semistable weights}

Consider a weight \(g:[-3,3] \to \bbR_{>0}\) which is even. Then it is obvious from the symmetry of the signed measure \(\nu\) under the reflection \(y\mapsto -y\) that \(\nu(g)=0\). 
If \(X\) were weight-insensitive K-polystable, then for any such weight, we would have \(\mu(g)>0\). 
Let us build an example where this is false. 

Let \(a\in \bbR\) be a parameter, and consider the family of weights \(g_a: y\mapsto \cosh(ay)\) indexed by \(a\). 
It is immediate that the function \(a\mapsto \mu(g_a)\) is continuous, and we have \(\mu(g_0)>0\) since \(Q^3\) admits a Kähler-Einstein metric. 
Actually, one may compute an exact expression:
\[ \mu(g_a) = \frac{-16}{a^4}\left(6a^2+a\sinh(3a)-2\cosh(3a)+2 \right) \]
and the fastest growing summand yields an equivalent as \(a\to +\infty\):
\[ \mu(g_a) \sim -\frac{16 e^{3a}}{a^3} \]
As a consequence, \(\mu(g_a)\) is negative for large values of \(a\). 
From intermediate value theorem, there is a value \(a_0\) such that \(\mu(g_{a_0})=0\) and we can give an approximate value: \(a_0\simeq 1.81037\). 

We have shown that for the weight \(g_a\), \(X\) is strictly weighted K-unstable if \(a\) is large enough, \(X\) is strictly weighted K-semistable if \(a=a_0\), and \(X\) is weighted K-polystable for values of \(a\) close to \(0\). 

\subsubsection{The toric degeneration}

By Theorem~\ref{thm_degenerations}, and the fact that the valuation cone of \(X\) is a half plane, there is a unique \(G\)-equivariant special test configuration for \(X\) up to twist. 
We denote its central fiber by \(Y\). 
As explained in \cite{Delcroix_2020} \(Y\) has the same moment polytope \(\Delta^+\) as \(X\) under the action of \(G\), its valuation cone is the full plane \(M\otimes \bbR\), and its weight lattice is still \(M\). 
By Theorem~\ref{thm_weighted_K-stab_spherical}, we obtain that, if \(X\) is strictly \(g\)-K-semistable for a weight \(g\), then \(Y\) is \(g\)-K-polystable. 
By uniqueness of such a degeneration, we have identified the optimal degeneration of \(X\) with respect to the weight \(g_{a_0}\) from the previous subsection. 
Note that since the optimal degeneration is obtained in a single step from a K-semistable to a K-polystable, there is no need for the assumption of log convexity in Linsheng Wang's paper \cite{Wang_optimal}. 

Let us now say a bit more about the variety \(Y\). 
A lot can be said from the theory of spherical varieties, we explain a bit here. 
By Pasquier \cite{Pasquier_2008}, since we know the moment polytope \(\Delta^+\) and weight lattice \(M\) with respect to the action of \(G\), we can fully identify \(Y\) as a rank two \(G\)-horospherical variety. 
Furthermore, \(Y\) admits an almost-faithful action of \(\SL_2\times \bbG_m^2\) via the action of \((G\times \Aut_G^0(Y))/Z(G)\), so in particular it admits an almost-faithful action of a three-dimensional torus. 
We deduce directly from this that \(Y\) is a toric \(\bbQ\)-Fano variety with degree 54. 
It follows rather straighforwardly from Brion's criterion for a divisor on a spherical variety to be Cartier \cite{Brion_1989} that \(Y\) is Gorenstein. 
By Pasquier \cite{Pasquier_2017}, we can also deduce for example that \(Y\) does not have terminal singularities. 
Additionally, from the horospherical description, we can readily see that the Fano index of \(Y\) is three, equal to its dimension. It then follows from \cite[Theorem~1.1]{Araujo-Druel_2014} that \(Y\) is a singular quadric which is consistent. 
From the classification of Gorenstein Fano threefolds \cite{Kreuzer-Skarke_1997, Kreuzer-Skarke_1998} as encoded in the Graded Ring Database \cite{grdb}, there would be a single candidate for \(Y\): the Gorenstein toric variety with reflexive ID 2. 

In a more down to earth approach, it suffices to exhibit an explicit non-product \(G\)-equivariant special test configuration for \(X\), to obtain an explicit description of \(Y\). 
As explained in the introduction, this explicit test configuration is the hypersurface \(\mathcal{X}\subset \bbP^4\times \bbC\) defined by the equation 
\[ x_0x_2-x_1^2 + zx_3x_4 = 0 \]
where \([x_0:\cdots:x_4]\) denote homogeneous coordinates on \(\bbP^4\), and \(z\) denotes the variable in \(\bbC\). 
The central fiber is thus the singular quadric defined by the equation 
\( x_0x_2-x_1^2=0\) 
in \(\bbP^4\), which indeed is the Gorenstein toric variety with reflexive ID 2, with respect to the action of \((\bbC^*)^3\) defined by 
\[ (t_1,t_2,t_3) \cdot [x_0:\cdots :x_4] = [t_1^2x_0:t_1x_1:x_2:t_2x_3:t_3x_4] \]

\subsection{The Fano threefold 2-19 and its toric optimal degeneration}

We now focus on the spherical action 3-2-19 on the Fano threefold 2-29. 
In this section \(X\) denotes the Fano threefold 2-29. 

\subsubsection{Unstable and semistable weights}

Consider a weight \(g:[-3,3] \to \bbR_{>0}\) which is even. Again, it is obvious from the symmetry of the signed measure \(\nu\) under the reflection \(y\mapsto -y\) that \(\nu(g)=0\). 
If 2-29 were weight-insensitive K-polystable, then for any such weight, we would have \(\mu(g)>0\). 
Let us build an example where this is false. 

Let \(a\in \bbR\) be a parameter, and consider again the family of weights \(g_a: y\mapsto \cosh(ay)\) indexed by \(a\). 
It is immediate that the function \(a\mapsto \mu(g_a)\) is continuous, and we have \(\mu(g_0)>0\) since 2-29 admits a Kähler-Einstein metric. 
The exact expression is now:
\[ \mu(g_a) = \frac{32}{a^4}\left(-\cosh(a)+\cosh(3a)-2a\sinh(a)-a\sinh(a)\cosh(2a)-a^2\cosh(a) \right) \]
and the fastest growing summand yields an equivalent as \(a\to +\infty\):
\[ \mu(g_a) \sim -\frac{32 e^{3a}}{a^3} \]
As a consequence, \(\mu(g_a)\) is negative for large values of \(a\). 
From intermediate value theorem, there is a value \(a_0\) such that \(\mu(g_{a_0})=0\) and we can give an approximate value: \(a_0\simeq 1.3176\). 

We have shown that for the weight \(g_a\), \(X\) is strictly weighted K-unstable if \(a\) is large enough, \(X\) is strictly weighted K-semistable if \(a=a_0\), and \(X\) is weighted K-polystable for values of \(a\) close to \(0\). 

\subsubsection{The toric degeneration}

The same explanations as for the quadric \(Q^3\) show that there is a unique \(G=\SL_2\times \bbG_m\)-equivariant special degeneration of \(X\) up to twist. 
Again, its central fiber \(Y\) is a Gorenstein, non-terminal, toric, \(g_{a_0}\)-K-polystable variety. 
Again, to identify \(Y\), it is enough to exhibit an explicit non-product special equivariant test configuration. 
For this we start from the test configuration \(\mathcal{X}\) for \(Q^3\) defined above. 
We blow up the linear section \(\mathcal{Z} \subset \mathcal{X}\) defined by \(x_3=x_4=0\), whose intersection with each fiber is the non-singular one dimensional \(G\)-stable sub-quadric. 
The resulting blowup \(\tilde{\mathcal{X}}\) is a non-product special test configuration for \(X\), whose central fiber is the blowup of the singular quadric above along the non-singular one-dimensional \(G\)-stable sub-quadric. 
From the toric point of view, we identify this central fiber as the Gorenstein toric variety with reflexive ID 19 in the Graded Ring Database \cite{grdb}. 

We have finished the proof of Theorem~\ref{thm:weight-sensitive}. 

\subsection{A strictly K-semistable log Fano pair which is weighted K-polystable}

We now consider the log Fano pairs \((X,tE)\) where \(X\) is the Fano manifold 2-29, and \(E\) is the exceptional divisor with respect to the blowup map \(X\to Q^3\). 
With the same conventions as in the previous sections, the moment polytope of this pair is the polytope with set of vertices \(\{\pm(0,3),\pm(4-2t,1+t)\}\). 
For example, for \(t=1/2\), this is: 

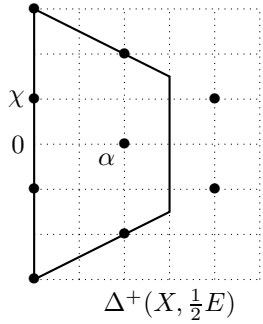
\begin{center}
\begin{tikzpicture}[scale=0.6]
\draw[dotted] (0,3) grid (5,-3);
\draw (0,0) node[left]{0};
\draw (0,1) node[left]{\(\chi\)};
\draw[thick] (0,3) -- (3,3/2) -- (3,-3/2) -- (0,-3) -- cycle;
\draw (2,0) node[below left]{\(\alpha\)};
\draw (0,3) node{\(\bullet\)}; \draw (0,1) node{\(\bullet\)}; \draw (0,-1) node{\(\bullet\)}; \draw (0,-3) node{\(\bullet\)}; \draw (2,0) node{\(\bullet\)}; \draw (2,-2) node{\(\bullet\)}; \draw (2,2) node{\(\bullet\)}; \draw (4,1) node{\(\bullet\)}; \draw (4,-1) node{\(\bullet\)}; 
\draw (3,-3) node[below]{\(\Delta^+(X,\frac{1}{2}E)\)};
\end{tikzpicture}
\end{center}

As before, \(g\)-weighted K-stability of \((X,tE)\) is encoded by two signed measures \(\nu_t\) and \(\mu_t\), with the conditions \(\nu_t(g)=0\) and \(\mu_t(g)>0\) corresponding to \(g\)-K-polystability. 
From symmetry, we have \(\nu_t(g)=0\) as soon as \(g\) is even, as for the manifold \(X\) itself. 
For even \(g\), we have  
\[ \mu_t(g) = \frac{8}{3} \left( (2-t)^2(1-2t)\int_0^{1+t}g(y)dy +\int_{1+t}^3g(y)(3-y)^2(3-2y)dy \right) \]
so in particular, 
\[ \mu_t(1) = -\frac{4}{3}(t-2)^2(3t^2+4t-2) \]
We deduce that, for \(t=t_0:=\frac{\sqrt{10}-2}{3}\), we have \(\mu_{t_0}(1)=0\). 
In other words, the log Fano pair \((X,t_0E)\) is strictly K-semistable. 

The optimal degeneration of this log Fano pair (in the classical sense), is by the same arguments as before the pair \((Y,tD)\), where \(Y\) is the Gorenstein toric threefold with reflexive ID 19, and \(D\) is the exceptional divisor in the blowup map from \(Y\) to the singular quadric which is the Gorenstein toric threefold with reflexive ID 2. 

Consider now the weight \(g\) defined by \(g(y)=\frac{1}{\cosh(y)}\). 
From numerical computations, one can easily check that \( \mu_{t_0}(g)>0 \), so that the pair \((X,t_0E)\) is weighted K-polystable with respect to this weight. 

Alternatively, to get rid of numerical computations, observe that the signed measure \(\mu_{t_0}\) is positive around zero. 
An obvious even function \(g\) such that \(\mu_{t_0}(g)>0\) is given by \(g=\indicator{[-\epsilon,\epsilon]}\) where \(\epsilon >0\) is small enough. 
Now observe that this non-smooth function may be arbitrarily well approximated by a smooth, even, positive, log concave function, this yields other examples of log concave weights for which the pair \((X,t_0E)\) is weighted K-polystable. 
We have finished the proof of Theorem~~\ref{thm:logpairs}.

\section{Higher dimensional examples}

We now quickly consider some higher dimensional examples, which are still spherical varieties of rank two, with a half plane as valuation cone. 
Let \(X\) be such a \(G\)-variety. 
From the same considerations as in the previous section, \(X\) admits a unique non-product \(G\)-equivariant special test configuration up to twist. 
The central fiber \(Y\) of this test configuration is a rank two horospherical variety. 
It will no longer, in general, be a toric variety. 
If \(X\) is strictly weighted K-semistable for a given weight \(g\), then \(Y\) will be weighted K-polystable for this weight, thus it will be the weighted optimal degeneration of \(X\). 

We consider the higher dimensional quadric \(Q^{n-2}\) with respect to the action of \(\SO_{n-2}(\bbC)\times \SO_2(\bbC)\subset \SO_n(\bbC)\), as studied in \cite[Section~4]{Delcroix_2022}. The combinatorial data needed to carry out the computations is available in this paper, we recall it in the following, adjusting the coordinates slightly so that the choices are consistent with those in the previous section (when \(n=5\)).  
In these coordinates \((x,y)\) (the \(x\) coordinate is dilated by four with respect to the one used in \cite{Delcroix_2022}, and the \(y\) coordinate is dilated by two), the moment polytope \(\Delta^+\) for the quadric \(Q^{n-2}\) is the triangle defined by the inequalities 
\(x\geq 0\), \(x\leq 2n-4+2t\) and \(x\leq 2n-4-2t\).

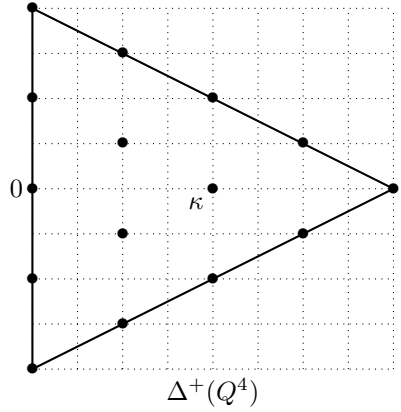
\begin{center}
\begin{tikzpicture}[scale=0.6]
\draw[dotted] (0,4) grid (8,-4);
\draw (0,0) node[left]{0};
\draw[thick] (0,4) -- (8,0) -- (0,-4) -- cycle;
\draw (4,0) node[below left]{\(\kappa\)};
\draw (0,4) node{\(\bullet\)}; \draw (0,2) node{\(\bullet\)}; \draw (0,0) node{\(\bullet\)}; \draw (0,-2) node{\(\bullet\)}; \draw (0,-4) node{\(\bullet\)}; \draw (2,1) node{\(\bullet\)}; \draw (2,-1) node{\(\bullet\)}; \draw (2,3) node{\(\bullet\)}; \draw (2,-3) node{\(\bullet\)}; \draw (4,0) node{\(\bullet\)}; \draw (4,-2) node{\(\bullet\)}; \draw (4,2) node{\(\bullet\)}; \draw (6,1) node{\(\bullet\)}; \draw (6,-1) node{\(\bullet\)}; \draw (8,0) node{\(\bullet\)};
\draw (4,-4) node[below]{\(\Delta^+(Q^{4})\)};
\end{tikzpicture}
\end{center}

The Duistermaat-Heckman polynomial is (a multiple of) \(x^{n-4}\), the point \(\kappa\) is the point with coordinates \((2n-8,0)\), and the valuation cone is the negative half plane defined by \(\kappa\) in \(N\otimes \bbR=(M\otimes \bbR)^*\). 
A weight for the action of \(\SO_2(\bbC)\simeq \bbC^*\) may, again, be considered as a function of the second variable \(y\). 
If such a weight \(g\) is even, then by symmetry the vanishing of the weighted Futaki invariant is automatic, and the weighted K-polystability condition is encoded by the positivity condition  
\begin{align*}
0 < & \int_{\Delta^+} (x-(2n-8))g(y)x^{n-4}dx dy \\
& = 2\int_0^{n-2} g(y) \left( \int_0^{2n-4-2y} (x-(2n-8))x^{n-4} dx\right) dy \\
& = \frac{4}{(n-2)(n-3)}\int_0^{n-2} g(y)(2n-4-2y)^{n-3}(n-2-(n-3)y) dy
\end{align*}
In the last expression, the affine function \(n-2-(n-3)y\) is positive for \(0\leq y < \frac{n-2}{n-3}\) and negative for \( \frac{n-2}{n-3} < y \leq n-2\). 
As a consequence, for any weight \(g\) that approximates well enough the function \(\indicator{[\frac{n-2}{n-3},n-2]}\), the quadric \(Q^{n-2}\) is weighted K-unstable. 
Since \(Q^{n-2}\) admits Kähler-Einstein metrics hence is K-polystable, we deduce that there exists weights for which \(Q^{n-2}\) is strictly weighted K-semistable, for which the optimal degeneration is given by the unique horospherical degeneration. 

Again, for \(Q^{n-2}\) it is easy to exhibit explicitly a non-product, \(\SO_{n-2}(\bbC)\times \SO_2(\bbC)\)-equivariant special test configuration, as the hypersurface defined by 
\[ x_0^2+\cdots +x_{n-3}^2 + z(x_{n-2}^2+x_{n-1}^2) =0 \]
in \(\bbP^{n-1}\times \bbC\). 
The central fiber is the singular quadric defined by the equation \(x_0^2+\cdots +x_{n-3}^2=0\) in \(\bbP^{n-1}\). 
We have thus proved Theorem~\ref{thm:high-dim}.

Considering the blowup of \(Q^{n-2}\) along the \(\SO_{n-2}(\bbC)\times \SO_2(\bbC)\)-stable subquadric \(Q^{n-4}\) of dimension \(n-4\) would again yield another Fano example which is not weight-insensitive K-polystable. We leave the details to the interested reader.

\bibliographystyle{alpha}
\bibliography{wsKs}

\end{document}